\documentclass{amsart}

\usepackage{amsmath,amsthm,amssymb,amsfonts,xcolor,url,kotex}
\usepackage[top=1in,bottom=1in,left=1.2in,right=1.2in]{geometry}

\usepackage[style=alphabetic]{biblatex}
\newcommand{\norm}[1]{\left\|#1\right\|}

\newtheorem{theorem}{Theorem}[subsection]
\renewcommand{\thetheorem}{%
    \ifnum\value{subsection}=0
        \thesection.\arabic{theorem}%
    \else
        \thesubsection.\arabic{theorem}%
    \fi
}

\usepackage[all]{xy}
\newcommand{\tp}{\mathop{\xymatrix{*+<.7ex>[o][F-]{\scriptstyle \top}}}}

\newtheorem{lemma}[theorem]{Lemma}

\theoremstyle{definition}

\theoremstyle{remark}
\newtheorem{remark}[theorem]{Remark}

\newtheorem*{maintheorem}{\normalfont\bfseries Main Theorem}

\title{A metric characterization of free unitary quantum groups}

\author{Yoonje Jeong}
\address{Yoonje Jeong\\
Department of Mathematics, University of California San Diego,
9500 Gilman Drive, La Jolla, CA 92093, USA}
\email{yojeong@ucsd.edu}

\author{Sang-Gyun Youn}
\address{Sang-Gyun Youn\\
Department of Mathematics Education and NextQuantum, Seoul National University,
Gwanak-Ro 1, Gwanak-Gu, Seoul, 08826, Republic of Korea}
\email{s.youn@snu.ac.kr}

\DeclareFieldFormat{postnote}{#1}
\DeclareCiteCommand{\cite}[\mkbibbrackets]
  {\usebibmacro{prenote}}
  {\usebibmacro{citeindex}\usebibmacro{cite}}
  {\multicitedelim}
  {\iffieldundef{postnote}
     {}
     {\setunit{\addcomma\space}\printfield{postnote}}}

\begin{document}

\begin{abstract}
For a compact matrix quantum group $G$ with a fundamental unitary representation $u$ of dimension at least two, the reduced operator norm of the self-adjoint main character satisfies $\norm{\chi_u+\chi_u^*}_{C_r(G)}\geq2\sqrt2$. Within the category of free quantum groups of Kac type, it is known that the minimum $2\sqrt2$ is attained only by free unitary quantum groups. We establish this characterization in full generality, including the non-Kac setting: this minimum is attained if and only if $G$ is isomorphic to a free unitary quantum group $U_Q^+$ as compact matrix quantum groups.
\end{abstract}

\maketitle

\section{Introduction}

Freeness is a fundamental concept that arises in many areas of mathematics. It naturally raises questions about what it means for one object to be `freer' than another, or for an object to be the `freest' within a given class. Group theory provides a standard algebraic formulation of these notions. Given two discrete groups $\Gamma_1$ and $\Gamma_2$, we regard $\Gamma_1$ as `freer' than $\Gamma_2$ if there exists a surjective group homomorphism $f:\Gamma_1\to\Gamma_2$. The free group $\mathbb{F}_n$ is then the `freest' among discrete groups generated by $n$ elements, since every such group is a quotient of $\mathbb{F}_n$. The operator-algebraic framework of locally compact quantum groups \cite{KV00,KV03} extends this perspective to finitely generated discrete quantum groups, where analogous notions of `freer' and `freest' arise. In this setting, free groups are no longer the universal objects, but the unitary free quantum groups $\mathbb{F}U_Q$ are.

Analytic approaches to these notions have been investigated in \cite{Ke59,CaCo22,JeYo26}. For finitely generated discrete quantum groups, the notion of `freer' is already well understood analytically through the moments of the self-adjoint main character \cite[Proposition 2.3]{JeYo26}. Characterizing the `freest' objects through operator norms, however, is more subtle. 

We begin by considering how the notion of `freest' can be understood analytically for finitely generated discrete groups. Let $\Gamma$ be a discrete group with a generating set $S=\{g_1,\ldots,g_n\}$, and let $\lambda:\Gamma\to\mathcal U(\ell^2(\Gamma))$ denote its left regular representation, given by
\begin{equation*}
    (\lambda(g)f)(x)=f(g^{-1}x),\qquad g,x\in\Gamma,\quad f\in\ell^2(\Gamma).
\end{equation*}
In this setting, the main character is $\displaystyle \chi_{\Gamma,S}=\sum_{i=1}^n\lambda(g_i)$. Then
\begin{equation*}
    2\sqrt{2n-1}\leq \norm{\chi_{\Gamma,S}+\chi_{\Gamma,S}^*}_{B(\ell^2(\Gamma))}\leq2n
\end{equation*}
for all $n\geq2$. Kesten's theorem \cite{Ke59} shows that the lower bound $2\sqrt{2n-1}$ is attained if and only if $g_1,\ldots,g_n$ freely generate $\Gamma$, i.e., $\Gamma\cong \mathbb{F}_n$. Thus, the `freest' groups are characterized by the minimum operator norm of the self-adjoint main character. A corresponding norm characterization of free Haar unitaries in finite von Neumann algebras was obtained by Cadilhac and Collins \cite{CaCo22}.

In the spirit of \cite{CaCo22,JeYo26}, let us describe the analogous minimization problem for quantum groups. There is a one-to-one correspondence between finitely generated discrete quantum groups and their dual compact matrix quantum groups, and this minimization problem can be formulated in the framework of compact matrix quantum groups. We denote by $U_Q^+$ the compact dual of the unitary free quantum group $\mathbb{F}U_Q$. Let $G$ be a compact matrix quantum group with a fundamental unitary representation $u=(u_{ij})_{1\leq i,j\leq n_u}$ of dimension $n_u\geq2$. The character of $u$ is $\displaystyle \chi_u=\sum_{i=1}^{n_u}u_{ii}$. By \cite[Proposition 3.1]{JeYo26}, we have
\[\norm{\chi_u+\chi_u^*}_{C_r(G)}\geq 2\sqrt{2}\]
for every such pair $(G,u)$, where $C_r(G)$ is the reduced $C^*$-algebra of $G$. The lower bound $2\sqrt{2}$ is attained by the free unitary quantum groups $U_Q^+$. The main question is whether these are the only minimizers. Within the category of free quantum groups of Kac type, a partial affirmative answer was obtained in \cite[Theorem 3.3]{JeYo26} using explicit descriptions of their representation categories.

In this paper, we resolve this question in full generality. We state our main result in the framework of compact matrix quantum groups.

\begin{maintheorem}
    Let $(G,u)$ be a compact matrix quantum group with fundamental unitary representation $u=(u_{ij})_{1\leq i,j\leq n_u}$, where $n_u\geq2$. Then
    \begin{equation*}
        \norm{\chi_u+\chi_u^*}_{C_r(G)}=2\sqrt2 \quad\Longleftrightarrow\quad (G,u)\cong(U_Q^+,v)
    \end{equation*}
    for some $Q\in\mathrm{GL}_{n_u}(\mathbb C)$. Here $v$ is the fundamental representation of $U_Q^+$, and the isomorphism identifies the specified fundamental matrix coefficients.
\end{maintheorem}

\subsection*{AI tools disclosure.}The main result of this paper was initially obtained by OpenAI's GPT-6 Astra. The authors subsequently verified, revised, and refined the proof and prepared the manuscript. The authors take full responsibility for the final text and all mathematical claims.


\section{Preliminaries}\label{sec2}

\subsection{Compact quantum groups and representations}

We recall the notation used below; for the basic representation theory, we refer to \cite{Woro87,NT13}. For a compact quantum group $G$, let $\operatorname{Pol}(G)$ denote its polynomial $*$-algebra, with comultiplication
\begin{equation*}
    \Delta_G:\operatorname{Pol}(G)\longrightarrow \operatorname{Pol}(G)\otimes\operatorname{Pol}(G).
\end{equation*}
The Haar state $h_G$ is characterized by
\begin{equation*}
    (h_G\otimes\mathrm{id})\Delta_G(a)=h_G(a)1 =(\mathrm{id}\otimes h_G)\Delta_G(a), \qquad a\in\operatorname{Pol}(G).
\end{equation*}
It is faithful on $\operatorname{Pol}(G)$, but need not be tracial.

A finite-dimensional unitary representation of $G$ on a Hilbert space $H_u$ of dimension $n_u$ is a unitary matrix $u=(u_{ij})_{1\leq i,j\leq n_u}\in B(H_u)\otimes\operatorname{Pol}(G)\cong M_{n_u}(\operatorname{Pol}(G))$ satisfying
\begin{equation*}
    \Delta_G(u_{ij})=\sum_{k=1}^{n_u}u_{ik}\otimes u_{kj}, \qquad 1\leq i,j\leq n_u.
\end{equation*}
Unless otherwise specified, all representations are finite-dimensional and unitary. The trivial representation is denoted by $1$, and the fixed-vector space of $u$ is
\begin{equation*}
    \operatorname{Fix}_G(u)=\{\xi\in H_u:u(\xi\otimes1)=\xi\otimes1\}.
\end{equation*}

For representations $u,v$ of $G$, their direct sum is the representation on $H_u\oplus H_v$ given by
\begin{equation*}
    u\oplus v=
    \begin{bmatrix}
        u&0\\
        0&v
    \end{bmatrix}.
\end{equation*}
Their tensor product, denoted by $u\tp v$, acts on $H_u\otimes H_v$ and has coefficients
\begin{equation*}
    (u\tp v)_{(i,k),(j,l)}=u_{ij}v_{kl}, \qquad 1\leq i,j\leq n_u,\quad 1\leq k,l\leq n_v.
\end{equation*}
The entrywise adjoint $u^c=(u_{ij}^*)_{1\leq i,j\leq n_u}$ need not be unitary, but is equivalent to a unitary representation. We write $\overline u$ for a unitary conjugate representation: choose $Q\in\mathrm{GL}_{n_u}(\mathbb C)$ such that
\begin{equation*}
    \overline u=Qu^cQ^{-1}
\end{equation*}
is unitary. Its equivalence class is independent of this choice. Every representation is a direct sum of irreducible representations. We denote their equivalence classes by $\operatorname{Irr}(G)$ and choose one representative of each class. Their matrix coefficients linearly span $\operatorname{Pol}(G)$.

A representation $u$ is called fundamental if its coefficients generate $\operatorname{Pol}(G)$ as a unital $*$-algebra. A compact quantum group equipped with such a representation is called a compact matrix quantum group and is denoted by $(G,u)$. For compact matrix quantum groups $(G,u)$ and $(H,v)$ with $n_u=n_v$, we say that $(G,u)$ is a matrix quantum subgroup of $(H,v)$ if there exists a surjective unital $*$-homomorphism $\pi:\operatorname{Pol}(H)\longrightarrow\operatorname{Pol}(G)$ satisfying
\begin{equation*}
    \pi(v_{ij})=u_{ij},\quad 1\leq i,j\leq n_v.
\end{equation*}
This map satisfies $\Delta_G\circ\pi=(\pi\otimes\pi)\circ\Delta_H$. We write $(G,u)\cong(H,v)$ when such a map is a $*$-isomorphism. Thus, throughout this paper, isomorphisms of compact matrix quantum groups identify the specified fundamental matrix coefficients and are understood at the level of polynomial $*$-algebras.

\subsection{Unitary free quantum groups}

Let $n\geq2$ and $Q\in\mathrm{GL}_{n}(\mathbb C)$. The unitary free quantum group $U_Q^+$ \cite{Wang95,VDW96} is the compact matrix quantum group with fundamental representation $v=(v_{ij})_{1\leq i,j\leq n}$ for which $\operatorname{Pol}(U_Q^+)$ is the universal unital $*$-algebra generated by the coefficients of $v$ subject to the relations that both $v$ and $Qv^cQ^{-1}$ are unitary. Its comultiplication is determined by
\begin{equation*}
    \Delta_{U_Q^+}(v_{ij})=\sum_{k=1}^{n}v_{ik}\otimes v_{kj},\qquad 1\leq i,j\leq n.
\end{equation*}

Given a finite-dimensional unitary representation $u$ of a compact quantum group $G$ with $n_u\geq2$, choose $Q\in\mathrm{GL}_{n_u}(\mathbb C)$ such that $Qu^cQ^{-1}$ is unitary. Then the universal property gives a unital $*$-homomorphism $\pi:\operatorname{Pol}(U_Q^+)\longrightarrow\operatorname{Pol}(G)$ satisfying
\begin{equation*}
\pi(v_{ij})=u_{ij},\quad 1\leq i,j\leq n_u.
\end{equation*}
It satisfies $\Delta_G\circ\pi =(\pi\otimes\pi)\circ\Delta_{U_Q^+}$. This map is surjective precisely when $u$ is fundamental. In that case, $(G,u)$ is a matrix quantum subgroup of $(U_Q^+,v)$, and $(G,u)\cong(U_Q^+,v)$ precisely when $\pi$ is injective.

Let $\mathcal W=\{\alpha,\beta\}^*$ be the free monoid on two letters, with empty word $e$. For $x\in\mathcal W$, write $|x|$ for its length. Define an involution by
\begin{equation*}
    \overline e=e,\qquad
    \overline\alpha=\beta,\qquad
    \overline\beta=\alpha,\qquad
    \overline{x_1\cdots x_k}
    =\overline{x_k}\cdots\overline{x_1}.
\end{equation*}
A word $x$ is called self-conjugate if $\overline x=x$. We use the following description of the irreducible representations of $U_Q^+$.

\begin{theorem}[{\cite[Theorem~1]{Banica97}}]\label{thm201}
    The irreducible representations of $U_Q^+$ are indexed by $\mathcal W$. Writing $r_x$ for the representation indexed by $x$, we have, up to equivalence,
    \begin{equation*}
        r_e=1,\qquad r_\alpha=v,\qquad
        r_\beta=\overline v,\qquad
        r_{\overline x}=\overline{r_x}.
    \end{equation*}
    For all $x,y\in\mathcal W$,
    \begin{equation}\label{eq201}
        r_x\tp r_y
        \cong\bigoplus_{\substack{a,b,g\in\mathcal W:\\
                    x=ag,\ y=\overline g b}}r_{ab}.
    \end{equation}
\end{theorem}

\begin{remark}\label{rem201}
    For nonempty $x,y\in\mathcal W$, a term with $g\neq e$ in \eqref{eq201} requires the first letter of $y$ to be the conjugate of the last letter of $x$. Thus, if the last letter of $x$ and the first letter of $y$ are equal, only the term $g=e$ occurs, and $r_x\tp r_y\cong r_{xy}$. We refer to this situation as having no cancellation at the junction of $x$ and $y$.
\end{remark}

\subsection{Characters and Haar integration}

The character of a representation $u$ of $G$ is
\begin{equation*}
    \chi_u=(\operatorname{Tr}\otimes\mathrm{id})(u)
    =\sum_{i=1}^{n_u}u_{ii},
\end{equation*}
where $\operatorname{Tr}$ is the unnormalized matrix trace. For representations $u,v$ of $G$,
\begin{equation*}
    \chi_{u\oplus v}=\chi_u+\chi_v,\qquad
    \chi_{u\text{\tiny $\tp$} v}=\chi_u\chi_v,\qquad
    \chi_{\overline u}=\chi_u^*.
\end{equation*}
In particular, sums, products, and adjoints of characters are again characters. We allow the zero character, corresponding to an empty direct sum. For $x\in\mathcal W$, write $\chi_x=\chi_{r_x}\in\operatorname{Pol}(U_Q^+)$. Thus $\chi_x$ always denotes a character on $U_Q^+$. The fusion rule gives
\begin{equation}\label{eq202}
    \chi_e=1,\qquad
    \chi_{\overline x}=\chi_x^*,\qquad
    \chi_x\chi_y
    =\sum_{\substack{a,b,g\in\mathcal W:\\
                    x=ag,\ y=\overline g b}}\chi_{ab}.
\end{equation}

Haar averaging gives the orthogonal projection $P_u=(\mathrm{id}\otimes h_G)(u)$ onto $\operatorname{Fix}_G(u)$. Consequently,
\begin{equation}\label{eq203}
    h_G(\chi_u)=\operatorname{Tr}(P_u)
    =\dim\operatorname{Fix}_G(u)\in\mathbb Z_{\geq0}.
\end{equation}
If $u$ is nontrivial and irreducible, then $P_u=0$, so $h_G(u_{ij})=0$ for all $1\leq i,j\leq n_u$.

By \cite[Theorem 5.8]{Woro87}, the irreducible characters satisfy
\begin{equation*}
    h_G(\chi_\gamma^*\chi_\eta)=\delta_{\gamma,\eta},
    \qquad \gamma,\eta\in\operatorname{Irr}(G).
\end{equation*}
Thus, if a nonzero representation $u$ has character
\begin{equation*}
    \chi_u=\sum_{\gamma\in\operatorname{Irr}(G)}m_\gamma\chi_\gamma,
\end{equation*}
with irreducible multiplicities $m_\gamma\in\mathbb Z_{\geq0}$, only finitely many of which are nonzero, then
\begin{equation*}
    h_G(\chi_u^*\chi_u)
    =\sum_{\gamma\in\operatorname{Irr}(G)}m_\gamma^2\geq1.
\end{equation*}
In particular,
\begin{equation}\label{eq204}
    \begin{aligned}
        h_G(\chi_u^*\chi_u)=1
        &\quad\Longleftrightarrow\quad u\text{ is irreducible},\\
        h_G(\chi_u^*\chi_u)\geq2
        &\quad\Longleftrightarrow\quad u\text{ is reducible}.
    \end{aligned}
\end{equation}
These relations hold without assuming that $h_G$ is tracial.

We shall also use
\begin{equation}\label{eq205}
    h_G(\chi_u\chi_v)\geq h_G(\chi_u)h_G(\chi_v).
\end{equation}
Indeed, write $\chi_u=h_G(\chi_u)1+\psi$ and $\chi_v=h_G(\chi_v)1+\eta$, where $\psi,\eta$ are characters satisfying $h_G(\psi)=h_G(\eta)=0$. Since $\psi\eta$ is also a character, \eqref{eq203} gives
\begin{equation*}
    h_G(\chi_u\chi_v)
    =h_G(\chi_u)h_G(\chi_v)+h_G(\psi\eta)
    \geq h_G(\chi_u)h_G(\chi_v).
\end{equation*}

Let $C_r(G)$ denote the completion of $\operatorname{Pol}(G)$ in the operator norm of the GNS representation of $h_G$. We identify $\operatorname{Pol}(G)$ with its image in $C_r(G)$. The Haar state extends to a faithful state on $C_r(G)$; in particular,
\begin{equation*}
    |h_G(a)|\leq\|a\|_{C_r(G)},\qquad a\in C_r(G).
\end{equation*}

\section{Proof of the main theorem}\label{sec3}

Let $G$ be a compact quantum group and let $u=(u_{ij})_{1\leq i,j\leq n_u}$ be a finite-dimensional unitary representation of $G$ with $n_u\geq2$. For the lower estimates below, $u$ need not be fundamental. Choose $Q\in\mathrm{GL}_{n_u}(\mathbb C)$ such that $Qu^cQ^{-1}$ is unitary. Then the universal property gives a unital $*$-homomorphism
\begin{equation*}
    \pi:\operatorname{Pol}(U_Q^+)\longrightarrow\operatorname{Pol}(G),
    \qquad \pi(v_{ij})=u_{ij},\quad 1\leq i,j\leq n_u.
\end{equation*}
This map $\pi$ is surjective when $u$ is fundamental. For $x\in\mathcal W$, write
\begin{equation*}
    r_x^\pi=(\mathrm{id}\otimes\pi)(r_x),\qquad
    \chi_x^\pi=\pi(\chi_x)=\chi_{r_x^\pi}\in\operatorname{Pol}(G).
\end{equation*}
In particular, $r_\alpha^\pi=u$ and $\chi_\alpha^\pi=\chi_u$. Although the $r_x$ are irreducible and pairwise inequivalent over $U_Q^+$, the representations $r_x^\pi$ need not be irreducible or pairwise inequivalent over $G$. Applying $\pi$ to \eqref{eq202}, we obtain
\begin{equation}\label{eq301}
    \chi_e^\pi=1,\qquad
    \chi_{\overline x}^\pi=(\chi_x^\pi)^*,\qquad
    \chi_x^\pi\chi_y^\pi
    =\sum_{\substack{a,b,g\in\mathcal W:\\
                    x=ag,\ y=\overline g b}}\chi_{ab}^\pi.
\end{equation}
For $k\geq0$, define
\begin{equation}\label{eq302}
    B_k=\sum_{\substack{x\in\mathcal W:\ |x|=k}}h_G(\chi_x^\pi).
\end{equation}
Every summand is a nonnegative integer by \eqref{eq203}, and $B_0=1$. If $\pi$ is injective, then $B_k=0$ for every $k\geq1$. We first establish a lower bound for $B_k$ when $\pi$ is not injective, beginning with the following lemma.

\begin{lemma}\label{lem301}
    If $\pi$ is not injective, then there exists a nonempty word $w_0\in\mathcal W$ such that
    \begin{equation*}
        \overline{w_0}=w_0,\qquad h_G(\chi_{w_0}^\pi)\geq1.
    \end{equation*}
\end{lemma}

\begin{proof}
    We first claim that $h_G(\chi_w^\pi)>0$ for some $w\in\mathcal W\setminus\{e\}$. Otherwise, for every $w\in\mathcal W\setminus\{e\}$, the orthogonal projection
    \begin{equation*}
        P_w^\pi=(\mathrm{id}\otimes h_G)(r_w^\pi)
        =(\mathrm{id}\otimes(h_G\circ\pi))(r_w)
    \end{equation*}
    satisfies
    \begin{equation*}
        \operatorname{Tr}(P_w^\pi)
        =(h_G\circ\pi)(\chi_w)=h_G(\chi_w^\pi)=0,
    \end{equation*}
    so $P_w^\pi=0$. Since $r_w$ is nontrivial and irreducible, we obtain
    \begin{equation*}
        (\mathrm{id}\otimes(h_G\circ\pi))(r_w)
        =P_w^{\pi}=0=(\mathrm{id}\otimes h_{U_Q^+})(r_w),
        \qquad w\in\mathcal W\setminus\{e\}.
    \end{equation*}
    Also, for $w=e$, we have
    \[(\mathrm{id}\otimes(h_G\circ\pi))(r_e)
        =1=(\mathrm{id}\otimes h_{U_Q^+})(r_e).\]
    Since the matrix coefficients of $r_w$'s span $\operatorname{Pol}(U_Q^+)$, it follows that $h_G\circ\pi=h_{U_Q^+}$ on $\operatorname{Pol}(U_Q^+)$. Thus, for every $a\in\ker\pi$,
    \begin{equation*}
        h_{U_Q^+}(a^*a)=h_G\bigl(\pi(a)^*\pi(a)\bigr)=0.
    \end{equation*}
    Then faithfulness of $h_{U_Q^+}$ on $\operatorname{Pol}(U_Q^+)$ gives $a=0$, contradicting the non-injectivity of $\pi$.

    Thus, we can take $w\in\mathcal W\setminus\{e\}$ such that $h_G(\chi_w^\pi)>0$. Since $w\neq e$, the classical dimension of $r_w$ is at least two by \cite[Lemmas~4.12 and~4.13]{KW22}. Thus $r_w^\pi$ has dimension at least two and a nonzero fixed vector, so it is reducible. By \eqref{eq204} and \eqref{eq301},
    \begin{equation*}
        2\leq h_G((\chi_w^\pi)^*\chi_w^\pi)
        =\sum_{\substack{a,g\in\mathcal W:\ w=ga}}
        h_G(\chi_{\overline a a}^\pi).
    \end{equation*}
    The summand with $a=e$ equals $h_G(1)=1$, and all summands are nonnegative. Hence $h_G(\chi_{\overline a a}^\pi)\geq1$ for some $a\in\mathcal W\setminus\{e\}$. Taking $w_0=\overline a a$ completes the proof.
\end{proof}

For the remainder of the lower-bound argument, assume that $\pi$ is not injective. Fix $w_0\in\mathcal W\setminus\{e\}$ as in Lemma~\ref{lem301}, and put
\begin{equation*}
    m_0=|w_0|,\qquad \lambda_0=h_G(\chi_{w_0}^\pi)\geq1.
\end{equation*}
Let $s_0$ and $t_0$ be the first and last letters of $w_0$, respectively. Since $\overline{w_0}=w_0$, we have $t_0=\overline{s_0}\neq s_0$. Define
\begin{equation*}
    \mathcal W_{t_0}
    =\{x\in\mathcal W:x\text{ begins with }t_0\}.
\end{equation*}

\begin{lemma}\label{lem302}
    For $x\in\mathcal W_{t_0}$, define $C(x)=\overline x w_0x$. Then $C(x)$ has length $m_0+2|x|$, is self-conjugate, and satisfies
    \begin{equation}\label{eq303}
        \chi_{C(x)}^\pi=(\chi_x^\pi)^*\chi_{w_0}^\pi\chi_x^\pi.
    \end{equation}
    Moreover,
    \begin{equation}\label{eq304}
        h_G(\chi_{C(x)}^\pi)
        \geq\lambda_0 h_G((\chi_x^\pi)^*\chi_x^\pi)
        \geq\lambda_0.
    \end{equation}
    If $x,y\in\mathcal W_{t_0}$ have different last letters, then
    \begin{equation}\label{eq305}
        \chi_{C(x)C(y)}^\pi
        =\chi_{C(x)}^\pi\chi_{C(y)}^\pi.
    \end{equation}
\end{lemma}

\begin{proof}
    By definition, $|C(x)|=m_0+2|x|$ and
    \begin{equation*}
        \overline{C(x)}
        =\overline x\,\overline{w_0}x
        =\overline x w_0x=C(x).
    \end{equation*}
    The last letter of $\overline x$ and the first letter of $w_0$ are both $s_0$, while the last letter of $w_0$ and the first letter of $x$ are both $t_0$. Thus there is no cancellation at either junction, and
    \begin{equation*}
        \chi_{C(x)}=\chi_{\overline x}\chi_{w_0}\chi_x=\chi_x^*\chi_{w_0}\chi_x.
    \end{equation*}
    Applying $\pi$, we obtain $\chi_{C(x)}^\pi=(\chi_x^\pi)^*\chi_{w_0}^\pi\chi_x^\pi$.

    Since $\lambda_0$ is the multiplicity of the trivial representation in $r_{w_0}^\pi$, write $\chi_{w_0}^\pi=\lambda_0 1+\psi$, where $\psi$ is a character with $h_G(\psi)=0$. The product $(\chi_x^\pi)^*\psi\chi_x^\pi$ is also a character, so its Haar integral is nonnegative. Therefore,
    \begin{align*}
        h_G(\chi_{C(x)}^\pi)
        &=\lambda_0 h_G((\chi_x^\pi)^*\chi_x^\pi)
          +h_G((\chi_x^\pi)^*\psi\chi_x^\pi)\\
        &\geq\lambda_0 h_G((\chi_x^\pi)^*\chi_x^\pi)
        \geq\lambda_0,
    \end{align*}
    where the last inequality follows from \eqref{eq204}.

    Finally, suppose that $x,y\in\mathcal W_{t_0}$ have different last letters. If $x$ ends in $b$, then $y$ ends in $\overline b$, so $C(x)$ ends in $b$ and $C(y)$ begins with $b$. Thus there is no cancellation at the junction of $C(x)$ and $C(y)$, and
    \begin{equation*}
        \chi_{C(x)C(y)}=\chi_{C(x)}\chi_{C(y)}.
    \end{equation*}
    Applying $\pi$, we obtain $\chi_{C(x)C(y)}^\pi
        =\chi_{C(x)}^\pi\chi_{C(y)}^\pi$.
\end{proof}

The following observation will be a key ingredient for Lemma \ref{lem304}. For a word $V=a_1\cdots a_L$, an integer $\delta$ with $1\leq\delta\leq L$ is called a \emph{period} of $V$ if $a_i=a_{i+\delta}$ whenever $1\leq i\leq L-\delta$.

\begin{lemma}\label{lem303}
    Suppose that $U,V\in\mathcal W$ satisfy $\overline U=U$, $\overline V=V$, and $V=UZ$ with $Z\neq e$. Then $|Z|$ is a period of $V$.
\end{lemma}

\begin{proof}
    Set $L=|V|$ and $\delta=|Z|$, so that $|U|=L-\delta$. Write $V=a_1\cdots a_L$ and $U=b_1\cdots b_{L-\delta}$. Comparing the first $|U|$ letters in $V=UZ$ gives
    \begin{equation*}
        a_i=b_i,\qquad 1\leq i\leq L-\delta,
    \end{equation*}
    while comparing the last $|U|$ letters in $V=\overline V=\overline ZU$ gives
    \begin{equation*}
        b_i=a_{i+\delta},\qquad 1\leq i\leq L-\delta.
    \end{equation*}
    Hence $a_i=a_{i+\delta}$ for every $1\leq i\leq L-\delta$, so $\delta$ is a period of $V$.
\end{proof}

For $\ell\geq2$ and $b\in\{\alpha,\beta\}=\{s_0,t_0\}$, put
\begin{equation*}
    \mathcal W_{t_0}(\ell,b)
    =\{x\in\mathcal W_{t_0}:|x|=\ell,\ x\text{ ends in }b\}.
\end{equation*}
Since the first and last letters are prescribed, $|\mathcal W_{t_0}(\ell,b)|=2^{\ell-2}$. For each integer $n\geq3$ and integers $0\leq p_1,p_2\leq n^2$, set
\begin{equation}\label{eq306}
    l_1=n^3+p_1,\qquad
    l_2=n^3+p_2,\qquad
    l_3=n^3-p_1-p_2.
\end{equation}
Then $l_1,l_2\geq n^3$ and $l_3\geq n^2(n-2)\geq9$, so all three lengths are at least nine, and $l_1+l_2+l_3=3n^3$. Let $\Omega_n$ be the set of all tuples $(p_1,p_2,x_1,x_2,x_3)$ with $0\leq p_1,p_2\leq n^2$ and
\begin{equation*}
    x_1\in\mathcal W_{t_0}(l_1,t_0),\qquad
    x_2\in\mathcal W_{t_0}(l_2,s_0),\qquad
    x_3\in\mathcal W_{t_0}(l_3,t_0).
\end{equation*}
Then $|\Omega_n|=(n^2+1)^2 2^{3n^3-6}$. Define $F_n:\Omega_n\to\mathcal W$ by
\begin{equation*}
    F_n(p_1,p_2,x_1,x_2,x_3)=C(x_1)C(x_2)C(x_3).
\end{equation*}
Every word in $\operatorname{Im}F_n$ has length $3m_0+6n^3$. Thus,
\begin{equation*}
    B_{3m_0+6n^3}
    =\sum_{\substack{z\in\mathcal W:\ |z|=3m_0+6n^3}}
      h_G(\chi_z^\pi)
    \geq\sum_{z\in\operatorname{Im}F_n}h_G(\chi_z^\pi).
\end{equation*}
For each $z\in\operatorname{Im}F_n$, choose $(p_1,p_2,x_1,x_2,x_3)\in\Omega_n$ such that $z=C(x_1)C(x_2)C(x_3)$. There is no cancellation at either junction, so Lemma~\ref{lem302} and \eqref{eq205} give
\begin{align*}
    h_G(\chi_z^\pi)=h_G(\chi_{C(x_1)}^\pi\chi_{C(x_2)}^\pi
          \chi_{C(x_3)}^\pi)\geq\prod_{j=1}^3 h_G(\chi_{C(x_j)}^\pi)
    \geq\lambda_0^3.
\end{align*}
Consequently,
\begin{equation}\label{eq307}
    B_{3m_0+6n^3}\geq\lambda_0^3|\operatorname{Im}F_n|.
\end{equation}
The following lemma estimates $|\operatorname{Im}F_n|$, and this plays a key role in establishing a lower bound of $B_k$.

\begin{lemma}\label{lem304}
    We have $\displaystyle \lim_{n\to\infty} \frac{|\operatorname{Im}F_n|}{|\Omega_n|}=1$.
\end{lemma}

\begin{proof}
    Define
    \begin{equation*}
        E_n=\{(\omega,\omega')\in\Omega_n^2:
        \omega\neq\omega',\ F_n(\omega)=F_n(\omega')\}.
    \end{equation*}
    For each $z\in\operatorname{Im}F_n$, put $\nu_z=|F_n^{-1}(z)|$. Then
    \begin{equation*}
        |\Omega_n|-|\operatorname{Im}F_n|
        =\sum_{z\in\operatorname{Im}F_n}(\nu_z-1)
        \leq\sum_{z\in\operatorname{Im}F_n}\nu_z(\nu_z-1)
        =|E_n|.
    \end{equation*}
    Thus it is enough to prove
    \begin{equation*}
        \lim_{n\to\infty}\frac{|E_n|}{|\Omega_n|}=0.
    \end{equation*}

    We now estimate $|E_n|$. For fixed integers $0\leq p_1,p_2,p_1',p_2'\leq n^2$, let $A(n,p_1,p_2,p_1',p_2')$ be the set of triples $(y_1,y_2,y_3)$ such that $\omega'=(p_1',p_2',y_1,y_2,y_3)\in\Omega_n$ and there exists $\omega=(p_1,p_2,x_1,x_2,x_3)\in\Omega_n$ satisfying
    \begin{equation}\label{eq308}
        C(x_1)C(x_2)C(x_3)
        =C(y_1)C(y_2)C(y_3)
    \end{equation}
    and $\omega\neq\omega'$.

    For the fixed parameters, the lengths $|C(x_j)|=m_0+2l_j$ are also fixed. Hence, by \eqref{eq308}, each $(y_1,y_2,y_3)$ determines $C(x_j)=\overline{x_j}w_0x_j$, and therefore $x_j$ for all $j=1,2,3$. Consequently,
    \begin{equation*}
        |E_n|
        =\sum_{0\leq p_1,p_2,p_1',p_2'\leq n^2}
        |A(n,p_1,p_2,p_1',p_2')|.
    \end{equation*}
    The same uniqueness holds with the two tuples interchanged, so we obtain
    \begin{equation*}
        |A(n,p_1,p_2,p_1',p_2')|
        =|A(n,p_1',p_2',p_1,p_2)|.
    \end{equation*}
    Note also that
    \begin{equation*}
        |y_1|=n^3+p_1',\qquad
        |y_2|=n^3+p_2',\qquad
        |y_3|=n^3-p_1'-p_2'.
    \end{equation*}
    We divide the argument into three cases according to whether $p_1=p_1'$ and, if so, whether $p_2=p_2'$.

    \smallskip
    \noindent Case 1: $p_1\neq p_1'$. By the symmetry above, it is enough to consider $p_1<p_1'$. In this case,
    \begin{equation*}
        |C(y_1)|-|C(x_1)|=2(p_1'-p_1)>0.
    \end{equation*}
    Comparing the first $|C(x_1)|$ letters on the two sides of \eqref{eq308}, we obtain
    \begin{equation*}
        C(y_1)=C(x_1)Z,\qquad |Z|=2(p_1'-p_1)>0
    \end{equation*}
    for some word $Z$. Since $\overline{C(x_1)}=C(x_1)$ and $\overline{C(y_1)}=C(y_1)$, Lemma~\ref{lem303} implies that $C(y_1)$ has period
    \begin{equation*}
        \delta=2(p_1'-p_1)\leq2n^2.
    \end{equation*}
    A word with period $\delta$ is determined by its first $\delta$ letters. Since both $|C(y_1)|$ and $\delta$ are fixed, there are at most $2^\delta$ choices for $C(y_1)$, and hence for $y_1$. Moreover,
    \begin{equation*}
        |y_2|+|y_3|
        =(n^3+p_2')+(n^3-p_1'-p_2')
        =2n^3-p_1'\leq2n^3.
    \end{equation*}
    Thus there are at most $2^{2n^3}$ choices for $(y_2,y_3)$, and we obtain
    \begin{equation*}
        |A(n,p_1,p_2,p_1',p_2')|
        \leq2^\delta 2^{2n^3}
        \leq2^{2n^3+2n^2}.
    \end{equation*}
    The same bound holds when $p_1>p_1'$ by interchanging the two tuples.

    \smallskip
    \noindent Case 2: $p_1=p_1'$ and $p_2\neq p_2'$. Since $|C(x_1)|=|C(y_1)|$, equation \eqref{eq308} gives
    \begin{equation*}
        C(x_1)=C(y_1),\qquad x_1=y_1,
    \end{equation*}
    and hence
    \begin{equation*}
        C(x_2)C(x_3)=C(y_2)C(y_3).
    \end{equation*}
    Suppose first that $p_2<p_2'$. Since
    \begin{equation*}
        |C(y_2)|-|C(x_2)|=2(p_2'-p_2)>0,
    \end{equation*}
    comparing the first $|C(x_2)|$ letters gives
    \begin{equation*}
        C(y_2)=C(x_2)Z,\qquad |Z|=2(p_2'-p_2)>0.
    \end{equation*}
    As in Case~1, Lemma~\ref{lem303} implies that $C(y_2)$ has period
    \begin{equation*}
        \delta=2(p_2'-p_2)\leq2n^2.
    \end{equation*}
    Thus there are at most $2^\delta$ choices for $y_2$. Moreover,
    \begin{equation*}
        |y_1|+|y_3|
        =(n^3+p_1)+(n^3-p_1-p_2')
        =2n^3-p_2'\leq2n^3.
    \end{equation*}
    Therefore,
    \begin{equation*}
        |A(n,p_1,p_2,p_1',p_2')|
        \leq2^\delta 2^{2n^3}
        \leq2^{2n^3+2n^2}.
    \end{equation*}
    The same bound holds when $p_2>p_2'$ by interchanging the two tuples.

    \smallskip
    \noindent Case 3: $p_1=p_1'$ and $p_2=p_2'$. In this case,
    \begin{equation*}
        |C(x_j)|=|C(y_j)|,\qquad j=1,2,3.
    \end{equation*}
    Comparing successive blocks in \eqref{eq308}, we obtain
    \begin{equation*}
        C(x_j)=C(y_j),\qquad x_j=y_j,\qquad j=1,2,3.
    \end{equation*}
    Thus $\omega=\omega'$, which is excluded by the definition of $A(n,p_1,p_2,p_1',p_2')$. Hence
    \begin{equation*}
        A(n,p_1,p_2,p_1',p_2')=\varnothing.
    \end{equation*}

    \smallskip
    There are exactly $(n^2+1)^4-(n^2+1)^2$ parameter quadruples satisfying $(p_1,p_2)\neq(p_1',p_2')$. Combining the three cases, we obtain
    \begin{equation}\label{eq309}
        |E_n|
        \leq\bigl((n^2+1)^4-(n^2+1)^2\bigr)
        2^{2n^3+2n^2}.
    \end{equation}
    Consequently,
    \begin{equation*}
        \frac{|E_n|}{|\Omega_n|}
        \leq\frac{\bigl((n^2+1)^4-(n^2+1)^2\bigr)
                     2^{2n^3+2n^2}}
                    {(n^2+1)^2 2^{3n^3-6}}
        =\frac{64n^2(n^2+2)}{2^{n^3-2n^2}}
        \longrightarrow0.
    \end{equation*}
    This proves the assertion.
\end{proof}

We now combine \eqref{eq307} with Lemma~\ref{lem304} to obtain a general lower bound for $B_k$ when $\pi$ is not injective.

\begin{theorem}\label{thm301}
    Suppose that $\pi$ is not injective. Then there exists $w_0\in\mathcal W\setminus\{e\}$ such that $\overline{w_0}=w_0$ and $\lambda_0=h_G(\chi_{w_0}^\pi)\geq1$. With $m_0=|w_0|$ and $c=\lambda_0^3/128$, we have
    \begin{equation*}
        B_{3m_0+6n^3}
        \geq\frac{\lambda_0^3}{128}(n^2+1)^2 2^{3n^3}
        \geq cn^4 2^{3n^3}
    \end{equation*}
    for all sufficiently large $n\in\mathbb N$. In particular,
    \begin{equation}\label{eq310}
        \limsup_{k\to\infty}
        \frac{B_k}{(k+1)2^{k/2}}=\infty.
    \end{equation}
\end{theorem}

\begin{proof}
    Fix $w_0$ as in Lemma~\ref{lem301} and apply the construction above. By \eqref{eq307},
    \begin{equation*}
        B_{3m_0+6n^3}\geq\lambda_0^3|\operatorname{Im}F_n|.
    \end{equation*}
    Lemma~\ref{lem304} gives
    \begin{equation*}
        |\operatorname{Im}F_n|
        \geq\frac12|\Omega_n|
        =\frac{(n^2+1)^2 2^{3n^3}}{128}
    \end{equation*}
    for all sufficiently large $n$. Therefore,
    \begin{equation*}
        B_{3m_0+6n^3}
        \geq\frac{\lambda_0^3}{128}(n^2+1)^2 2^{3n^3}
        \geq cn^4 2^{3n^3}.
    \end{equation*}
    Since $m_0$ and $c>0$ are independent of $n$,
    \begin{equation*}
        \frac{B_{3m_0+6n^3}}
        {(3m_0+6n^3+1)2^{(3m_0+6n^3)/2}}
        \geq\frac{cn^4}{(3m_0+6n^3+1)2^{3m_0/2}}
        \longrightarrow\infty.
    \end{equation*}
    This proves \eqref{eq310}.
\end{proof}

We now derive an upper bound for $B_k$ under the norm condition $\|\chi_u+\chi_u^*\|_{C_r(G)}=2\sqrt2$ and combine it with Theorem~\ref{thm301} to prove the main theorem.

\begin{theorem}\label{thm302}
    Let $(G,u)$ be a compact matrix quantum group with fundamental unitary representation $u=(u_{ij})_{1\leq i,j\leq n_u}$, where $n_u\geq2$. Choose $Q\in\mathrm{GL}_{n_u}(\mathbb C)$ such that $Qu^cQ^{-1}$ is unitary, and let $\pi:\operatorname{Pol}(U_Q^+)\longrightarrow\operatorname{Pol}(G)$ be the unital $*$-homomorphism determined by $\pi(v_{ij})=u_{ij}$ for $1\leq i,j\leq n_u$, where $v$ is the fundamental representation of $U_Q^+$. If $\|\chi_u+\chi_u^*\|_{C_r(G)}=2\sqrt2$, then
    \begin{equation*}
        B_k=\sum_{\substack{x\in\mathcal W: |x|=k}}h_G(\chi_x^\pi)
        \leq(k+1)2^{k/2},\qquad k\geq0,
    \end{equation*}
    and $\pi$ is a $*$-isomorphism. Consequently, $(G,u)\cong(U_Q^+,v)$ as compact matrix quantum groups.
\end{theorem}

\begin{proof}
Put $X_u=\chi_u+\chi_u^*=\chi_\alpha^\pi+\chi_\beta^\pi$. By \eqref{eq301}, we have
    \begin{equation*}
        X_u\cdot \sum_{\substack{x\in\mathcal W: |x|=k}}\chi_x^\pi
        =\sum_{\substack{x\in\mathcal W: |x|=k+1}}\chi_x^\pi
        +2\sum_{\substack{x\in\mathcal W: |x|=k-1}}\chi_x^\pi,
        \qquad k\geq1.
    \end{equation*}
    Indeed, the terms without cancellation give each word of length $k+1$ once, while each word of length $k-1$ occurs once for each of the two possible cancelled letters. It follows by induction that
    \begin{equation*}
        \sum_{\substack{x\in\mathcal W:\\ |x|=k}}\chi_x^\pi
        =P_k(X_u),\qquad k\geq0,
    \end{equation*}
    where the polynomials $P_k$ are defined by $P_0(t)=1$, $P_1(t)=t$, and
    \begin{equation*}
        P_{k+1}(t)=tP_k(t)-2P_{k-1}(t),\quad k\geq1.
    \end{equation*}

    Let $U_k$ denote the Chebyshev polynomial of the second kind, characterized by
    \begin{equation*}
        U_k(\cos\theta)
        =\frac{\sin((k+1)\theta)}{\sin\theta},
        \qquad 0<\theta<\pi.
    \end{equation*}
    These polynomials satisfy $U_0(t)=1$, $U_1(t)=2t$, and
    \begin{equation*}
        U_{k+1}(t)=2tU_k(t)-U_{k-1}(t),\quad k\geq1.
    \end{equation*}
    Then comparing the recurrences gives
    \begin{equation*}
        P_k(t)=2^{k/2}U_k\left(\frac{t}{2\sqrt2}\right).
    \end{equation*}
    Consequently,
    \begin{equation*}
        \max_{|t|\leq2\sqrt2}|P_k(t)|
        =2^{k/2}\max_{|s|\leq1}|U_k(s)|
        =(k+1)2^{k/2}.
    \end{equation*}

    Since $X_u$ is self-adjoint and $\|X_u\|_{C_r(G)}=2\sqrt2$, we have
    \begin{equation*}
        \sigma_{C_r(G)}(X_u)\subseteq[-2\sqrt2,2\sqrt2].
    \end{equation*}
    Thus continuous functional calculus gives
    \begin{align*}
        B_k=h_G(P_k(X_u))&\leq\|P_k(X_u)\|_{C_r(G)}
         =\max_{t\in\sigma_{C_r(G)}(X_u)}|P_k(t)|\\
        &\leq\max_{|t|\leq2\sqrt2}|P_k(t)|
         =(k+1)2^{k/2}.
    \end{align*}
This implies $\displaystyle \limsup_{k\to\infty} \frac{B_k}{(k+1)2^{k/2}}\leq1$, so $\pi$ is injective by Theorem \ref{thm301}. Since $u$ is fundamental, $\pi$ is also surjective, and hence is a $*$-isomorphism.
\end{proof}

\begin{remark}
    One might wonder whether a similar argument applies when $\|\chi_u+\chi_u^*\|_{C_r(G)}=2\sqrt2+\epsilon$ for sufficiently small $\epsilon>0$. Let $\theta_\epsilon>0$ be determined by $2\sqrt2+\epsilon=2\sqrt2\cosh\theta_\epsilon$. For every $k\geq0$, we have
    \begin{align*}
        \max_{|t|\leq2\sqrt2+\epsilon}|P_k(t)| =2^{k/2}
          \frac{\sinh((k+1)\theta_\epsilon)}
               {\sinh\theta_\epsilon}=(\sqrt2 e^{\theta_\epsilon})^k
          \frac{1-e^{-2(k+1)\theta_\epsilon}}
               {1-e^{-2\theta_\epsilon}}
          \sim(\sqrt2 e^{\theta_\epsilon})^k.
    \end{align*}
The last estimate holds up to positive constants that may depend on $\epsilon$ but are independent of $k$. Thus, repeating the argument used in the proof of Theorem~\ref{thm302} does not yield a contradiction in this setting.
\end{remark}

\emph{Acknowledgements}: S.-G.Y. was supported by the National Research Foundation of Korea (NRF) grant funded by the Korea government (MSIT) (No. RS-2025-00561391, No. RS-2025-10802969 and No. RS-2024-00413957).

\printbibliography

\end{document}